\documentclass{llncs}
\usepackage{amsmath}
\usepackage{amssymb}

\newcommand{\one}{\mathbf{1}}
\newcommand{\zero}{\mathbf{0}}

\newcommand{\R}{\mathbb{R}}
\newcommand{\Z}{\mathbb{Z}}
\newcommand{\N}{\mathbb{N}}
\DeclareMathOperator{\expOp}{Exp}
\newcommand{\expect}[1]{\expOp[\,{#1}\,]}
\newcommand{\expectcond}[2]{\expOp[\,{#1}\ |\ {#2}\,]}
\DeclareMathOperator{\probOp}{Prob}
\newcommand{\prob}[1]{\probOp[\,{#1}\,]}
\newcommand{\probcond}[2]{\probOp[\,{#1}\ |\ {#2}\,]}
\DeclareMathOperator{\convOp}{conv}
\newcommand{\conv}[1]{\convOp({#1})}

\DeclareMathOperator{\distOp}{dist}
\newcommand{\dist}[2]{\distOp({#1},{#2})}
\DeclareMathOperator{\edgesOp}{E}
\newcommand{\edges}[1]{\edgesOp({#1})}
\DeclareMathOperator{\vertsOp}{V}
\newcommand{\verts}[1]{\vertsOp({#1})}

\DeclareMathOperator{\pmcubeOp}{Q}
\newcommand{\pmcube}[1]{\pmcubeOp_{#1}}
\DeclareMathOperator{\pmcubestarOp}{Q^{\star}}
\newcommand{\pmcubestar}[1]{\pmcubeOp^{\star}_{#1}}

\newcommand{\setdef}[2]{\{\ {#1}\ :\ {#2}\ \}}
\newcommand{\bigsetdef}[2]{\big\{\ {#1}\ :\ {#2}\ \big\}}
\newcommand{\narrowsetdef}[2]{\{{#1}:{#2}\}}

\newcommand{\vecconf}[1]{\mathcal{V}_{#1}}

\DeclareMathOperator{\ddOp}{d}
\newcommand{\dd}[1]{\ddOp\!{#1}}
\newcommand{\lowfac}[2]{({#1})_{#2}}
\DeclareMathOperator{\const}{const}
\newcommand{\Euler}{e}

\DeclareMathOperator{\cutOp}{CUT}
\newcommand{\cut}[1]{\cutOp_{#1}}

\newtheorem{observation}{Observation}
\begin{document}

%
\mainmatter              

\title{On the Graph-Density of Random 0/1-Polytopes}
\titlerunning{Random 0/1-Polytopes}

\author{Volker Kaibel\inst{1} \and Anja Remshagen\inst{2}}
\authorrunning{V. Kaibel and A. Remshagen}
\tocauthor{Volker Kaibel (TU Berlin), Anja Remshagen (Univ. of West Georgia)}

\institute{%
  DFG Research Center ``Mathematics for key technologies'',
  MA~6--2, 
  TU~Berlin, 
  Stra\ss e des 17.~Juni~136, 
  10623~Berlin, 
  Germany, 
  \email{kaibel@math.tu-berlin.de}
  \and
  Department of Computer Science, 
  State University of West Georgia, 
  1600~Maple~Street, 
  Carrollton, GA~30118, 
  USA, 
  \email{anja@westga.edu}
}
\maketitle

\begin{abstract}
  Let $X_{d,n}$ be an $n$-element subset of $\{0,1\}^d$ chosen
  uniformly at random, and denote by $P_{d,n}:=\convOp X_{d,n}$ its convex
  hull. Let $\Delta_{d,n}$ be the density of the graph of $P_{d,n}$
  (i.e., the number of one-dimensional faces of $P_{d,n}$ divided
  by~$\binom{n}{2}$). Our main result is that, for any function $n(d)$, the
  expected value of $\Delta_{d,n(d)}$ converges (with
  $d\rightarrow\infty$) to one if, for some arbitrary
  $\varepsilon>0$, $n(d)\le (\sqrt{2}-\varepsilon)^d$ holds for
  all large~$d$, while it converges to zero if $n(d)\ge
  (\sqrt{2}+\varepsilon)^d$ holds for all large~$d$.
\end{abstract}

\section{Introduction}
\label{sec:intro}

Polytopes whose vertices have coordinates in $\{0,1\}$
(\emph{0/1-polytopes}) are the objects of study in large parts of
polyhedral combinatorics (see~\cite{Sch03}). Since that theory has
started to grow, people have been interested in the \emph{graphs}
(defined by the vertices and the one-dimensional faces) of the
polytopes under investigation. The main reason for this interest was,
of course, the role played by polytope graphs with respect to linear
programming and, in particular, the simplex algorithm.

Later it was recognized that the graphs of the 0/1-polytopes
associated with certain combinatorial objects (such as  matchings in a
graph or bases of a matroid) might also yield good candidates for
neighborhood structures with respect to the construction of random
walks for random generation of the respective objects. A quite
important (yet unsolved) problem arising in this context is the
question whether the graphs of 0/1-polytopes have good expansion
properties (see \cite{FM92,JS97,Kai01}).

We are short of knowledge on the graphs of \emph{general}
0/1-polytopes~\cite{Zie00}.  Among the few
exceptions are results about their diameters \cite{Nad89} and their
cycle structures \cite{NP84}.  Particularly striking is the fact that
several \emph{special} 0/1-polytopes associated with combinatorial
problems have quite dense graphs. The most prominent example for this
is probably the \emph{cut polytope} $\cut{k}$, i.e., the convex hull
of the characteristic vectors of those subsets of edges of the
complete graph $K_k$ that form cuts in $K_k$. Barahona and
Mahjoub~\cite{BM86} proved that the graph of $\cut{k}$ is complete,
i.e., its density equals one (where the \emph{density} of a graph
$G=(V,E)$ is $|E|/\binom{|V|}{2}$).  Since the dimension of $\cut{k}$
is $d=\binom{k}{2}$ and there are $n=2^{k-1}$ cuts in $K_k$, the cut
polytopes yield an infinite series of $d$-dimensional 0/1-polytopes
with (roughly) $c^{\sqrt{d}}$ vertices (for some constant $c$) and
graph-density one.

In this paper, we investigate the question for the graph-density of a
typical (i.e., random) 0/1-polytope. The (perhaps surprising) result
is that in fact the high density of the graphs of several
0/1-polytopes important in polyhedral combinatorics (such as the cut
polytopes) is not atypical at all. Our main result is the following
theorem, where $\expect{\cdot}$ denotes the expected value.

\begin{theorem}
  \label{thm:dens}
  Let $n:\N\longrightarrow\N$ be a function, and let
  $P_{d,n(d)}:=\convOp X_{d,n(d)}$ with an $n(d)$-element subset
  $X_{d,n(d)}$ of $\{0,1\}^d$ that is chosen uniformly at random.
  Denote by $\Delta_{d,n(d)}$ the density of the graph of
  $P_{d,n(d)}$.
  \begin{enumerate}
  \item[(i)] If there is some $\varepsilon>0$ such that $n(d)\le
    (\sqrt{2}-\varepsilon)^d$ for all sufficiently large $d$, then    
    $
    \displaystyle\lim_{d\rightarrow\infty}\expect{\Delta_{d,n(d)}}=1
    $.
  \item[(ii)] If there is some $\varepsilon>0$ such that $n(d)\ge
    (\sqrt{2}+\varepsilon)^d$ for all sufficiently large $d$,
    then 
    $
    \displaystyle\lim_{d\rightarrow\infty}\expect{\Delta_{d,n(d)}}=0
    $.
  \end{enumerate}
\end{theorem}

There is a similar threshold phenomenon for the volumes of random
0/1-polytopes.  Let $\Tilde{P}_{d,n(d)}$ be the convex hull of $n(d)$
points in $\{0,1\}^d$ that are chosen independently uniformly at
random (possibly with repetitions).  Dyer, F\"uredi, and
McDiarmid~\cite{DFM92} proved that the limit (for
$d\rightarrow\infty$) of the expected value of the $d$-dimensional
volume of $\Tilde{P}_{d,n(d)}$ is zero if, for some $\varepsilon>0$,
$n(d)\le(\frac{2}{\sqrt{e}}-\varepsilon)^d$ holds for all 
sufficiently large $d$, and it is one if, for some $\varepsilon>0$,
$n(d)\ge(\frac{2}{\sqrt{e}}+\varepsilon)^d$ holds for all 
sufficiently large~$d$.
Due to $\frac{2}{\sqrt{e}}<1.214$ and $\sqrt{2}>1.414$, one can deduce
(we omit the details) the following result from this and
Theorem~\ref{thm:dens}. It may be a bit surprising due to the fact
that the only $d$-dimensional 0/1-polytope with $d$-dimensional volume
equal to one is the 0/1-cube $\convOp\{0,1\}^d$, which has
only graph-density $\frac{d}{2^d-1}$.

\begin{corollary}
  For every $\delta>0$ there are (infinitely many) 0/1-polytopes with
  both graph density and volume at least $(1-\delta)$.
\end{corollary}

Another threshold result that is related to our work is due to
F\"uredi~\cite{Fue86}. He showed that, in the setting of
Theorem~\ref{thm:dens}, the limit (for $d\rightarrow\infty$) of the
probability that $P_{d,n(d)}$ contains the center of the 0/1-cube is
zero if, for some $\varepsilon>0$, $n(d)\le(2-\varepsilon)\cdot d$
holds for all sufficiently large $d$, and it is one if, for some
$\varepsilon>0$, $n(d)\ge(2+\varepsilon)\cdot d$ holds for all
sufficiently large $d$.  The material in
Sections~\ref{subsec:hyperplane}, \ref{subsec:bounds},
and~\ref{subsec:thresh} of our paper is very much inspired by
F\"uredi's work.

The aim of Sections~\ref{sec:tau} and~\ref{sec:p}
is to prove Theorem~\ref{thm:dens}.
Since it is a bit more convenient, we switch
from 0/1-polytopes to polytopes whose vertices have coordinates in
$\{-1,+1\}$ (\emph{$\pm 1$-polytopes}). Recalling that the density of
a graph equals the probability of a randomly chosen pair of its nodes
to be adjacent, 
Propositions~\ref{prop:pthrlow} and~\ref{prop:pthrup}
(Section~\ref{sec:p}), together with Proposition~\ref{prop:pmonotone},
imply Theorem~\ref{thm:dens} (with the $\varepsilon$'s in
Propositions~\ref{prop:pthrlow} and~\ref{prop:pthrup} replaced by
$\log\frac{\sqrt{2}}{\sqrt{2}-\varepsilon}$ and
$\log\frac{\sqrt{2}+\varepsilon}{\sqrt{2}}$, respectively).


We close with a few remarks in Section~\ref{sec:rem}.

\section{The Long-Edge Probability \boldmath{$\tau(k,m)$}}
\label{sec:tau}

We define $\pmcube{d}:=\{-1,+1\}^d$ and
$\pmcubestar{d}:=\pmcube{d}\setminus\{-\one,\one\}$ (where $\one$ is
the all-one vector). For $v,w\in\pmcube{d}$, 
denote by $\pmcubeOp(v,w)$ the subset of all points in $\pmcube{d}$
that agree with~$v$ and~$w$ in all components, where~$v$ and~$w$
agree. Thus, $\pmcubeOp(v,w)$ is the vertex set of the smallest face
of $\convOp \pmcube{d}$ containing~$v$ and~$w$. The dimension of this
face is 
$$
\dist{v}{w}\ :=\ \#\setdef{i\in\{1,\dots,d\}}{v_i\not= w_i}
$$
(the \emph{Hamming distance} of~$v$ and~$w$). Let
$\pmcubestarOp(v,w):=\pmcubeOp(v,w)\setminus\{v,w\}$.

We refer to~\cite{Zie95} for all notions and results from polytope
theory that we rely on. For a polytope~$P$, we denote by $\verts{P}$
and $\edges{P}$ the sets of vertices and edges of~$P$, respectively.
Recall that, for $X\subseteq\pmcube{d}$, we have $\verts{\convOp
  X}=X$. 

The following fact is essential for our treatment. It can easily be
deduced from elementary properties of convex polytopes.

\begin{lemma}
\label{lem:adjacency}
  For two vertices~$v$ and~$w$ of a $\pm 1$-polytope~$P\subset\R^d$
  we have
  $$
    \{v,w\}\in\edges{P}
    \ \Longleftrightarrow\ 
    \convOp\{v,w\}\cap\conv{P\cap\pmcubestarOp(v,w)}=\varnothing\ .
  $$
\end{lemma}

Throughout this section, let $Y_{k,m}\in\binom{\pmcubestar{k}}{m}$
(the $m$-element subsets of $\pmcubestar{k}$) be drawn uniformly at
random and define
$$
\tau(k,m)\ :=\ \prob{\conv{Y_{k,m}}\cap\convOp\{-\one,\one\}=\varnothing}\ .
$$
Thus, $\tau(k,m)$ is the probability that the ``long edge''
$\convOp\{-\one,\one\}$ is an edge of the polytope
$\conv{Y_{k,m}\cup\{-\one,\one\}}$. 
The next lemma follows from Lemma~\ref{lem:adjacency}.

\begin{lemma}
  \label{lem:tau}
  Let $X_{d,n}\in\binom{\pmcube{d}}{n}$ be chosen uniformly at random,
  defining the polytope $P_{d,n}:=\convOp X_{d,n}$. Choose
  a two-element subset $\{v,w\}$ of $X_{d,n}$ uniformly at random. Then,
  for every $k\in\{1,\dots,d\}$ and $m\in\{0,\dots,\min\{2^k-2,n-2\}\}$, we have
  the equation
  $$
  \probcond{\{v,w\}\in\edges{P_{d,n}}}{\dist{v}{w}=k,\#(X_{d,n}\cap\pmcubestarOp(v,w))=m}\ =\ \tau(k,m)\ .
  $$
\end{lemma}
 
Via Lemma~\ref{lem:tau}, asymptotic bounds on $\tau(k,m)$ will turn
out to be important for the proofs in Section~\ref{sec:p}.  In fact,
we will basically compute (or estimate) the probability $\pi(d,n)$
(see Section~\ref{sec:p}) that two randomly chosen vertices of a
$d$-dimensional random $\pm 1$-polytope with~$n$ vertices are adjacent
by partitioning the probability space into the events 
``$\dist{v}{w}=k$ and $\#(X_{d,n}\cap\pmcubestarOp(v,w))=m$''
for all $k\in\{1,\dots,d\}$ and $m\in\{0,\dots,\min\{2^k-2,n-2\}\}$.

For the study of $\tau(k,m)$, it is convenient to consider the
conditional probability
$$
\alpha(k,m)\ :=\ \probcond{\conv{Y_{k,m}}\cap\convOp\{-\one,\one\}=\varnothing}{Y_{k,m}\cap(-Y_{k,m})=\varnothing}\ ,
$$
which is related to $\tau(k,m)$ in the following way.

\begin{lemma}
  \label{lem:tau2}
  For $0\le m\le 2^k-2$ we have
  $$
  \tau(k,m)\ =\ \frac{\binom{2^{k-1}-1}{m}\cdot 2^m}{\binom{2^k-2}{m}}\cdot\alpha(k,m)\ .
  $$
\end{lemma}

\begin{proof}
  Clearly, $\conv{Y_{k,m}}\cap\convOp\{-\one,\one\}=\varnothing$
  implies $Y_{k,m}\cap(-Y_{k,m})=\varnothing$. Thus, the statement in
  the lemma is due to the fact that the number of sets
  $Y'\in\binom{\pmcubestar{k}}{m}$ with $Y'\cap(-Y')=\varnothing$ is
  $\binom{2^{k-1}-1}{m}\cdot 2^m$.
\end{proof}

We will first show that $\alpha(k,m)$ can be interpreted as a
conditional probability that a random $m$-element subset of a certain
vector configuration in $\R^{k-1}$ does not contain the origin in its
convex hull (Section~\ref{subsec:vecconfig}).  The latter probability
is then related to the expected number of chambers in a certain random
hyperplane arrangement. This number of chambers is finally estimated
via a well-known bound due to Harding
(Section~\ref{subsec:hyperplane}).

As a point of reference for the proofs in Section~\ref{sec:p}, let us
state the following monotonicity result here, whose (straightforward)
proof we omit.

\begin{lemma}
  \label{lem:taumono}
  For $0\le m\le 2^k-3$, we have
  $
  \tau(k,m)\ge\tau(k,m+1)
  $.
\end{lemma}

\subsection{The Vector Configuration \boldmath{$\vecconf{r}$}}
\label{subsec:vecconfig}

Let $\varphi:\R^{r+1}\longrightarrow H_{\one}\longrightarrow\R^{r}$ 
denote the orthogonal projection of $\R^{r+1}$ onto the
hyperplane $H_{\one}:=\narrowsetdef{x\in\R^{r+1}}{\one^Tx=0}$, followed
by the orthogonal projection to the first $r$ coordinates. We denote
by $\vecconf{r}:=\varphi(\pmcubestar{r+1})$ the image of
$\pmcubestar{r+1}$ under the projection~$\varphi$. We omit the simple
proof of the following result.

\begin{lemma}
  \label{lem:varphi}
  The projection~$\varphi$ is one-to-one on $\pmcubestar{r+1}$.
\end{lemma}

\begin{lemma}
\label{lem:ttau}
  For $Z_{r,m}\in\binom{\vecconf{r}}{m}$ chosen uniformly at random, we have
  $$
  \alpha(r+1,m)\ =\ \probcond{\zero\not\in\conv{Z_{r,m}}}{Z_{r,m}\cap(-Z_{r,m})=\varnothing}\ .
  $$
\end{lemma}

\begin{proof}
  Since
  $
  \convOp Y_{k,m}\cap\convOp\{-\one,\one\}=\varnothing
  $
  holds if and only if
  $
  \zero\not\in\convOp\varphi(Y_{k,m})
  $ holds,
  the claim follows from Lemma~\ref{lem:varphi} (because
  $Y_{k,m}\cap(-Y_{k,m})=\varnothing$ is equivalent to
  $\varphi(Y_{k,m})\cap(-\varphi(Y_{k,m}))=\varnothing$).
\end{proof}

With 
$\vecconf{r}^+:=\varphi\narrowsetdef{v\in\pmcubestar{r+1}}{v_{r+1}=+1}$,
we have $\vecconf{r}=\vecconf{r}^+\cup(-\vecconf{r}^+)$ and
$\vecconf{r}^+\cap(-\vecconf{r}^+)=\varnothing$.
For any fixed finite subset $S\subset\R^r$, and 
a uniformly at random chosen $\varepsilon\in\{-1,+1\}^S$, denote
$
\alpha(S):=\prob{\zero\not\in\convOp\narrowsetdef{\varepsilon_s s}{s\in S}}
$.

\begin{lemma}
  \label{lem:alphazplus}
  Let $Z_{r,m}^+\in\binom{\vecconf{r}^+}{m}$ be chosen uniformly at
  random. Then we have
  $$
  \alpha(r+1,m)\ =\ \expect{\alpha(Z_{r,m}^+)}\ .
  $$
\end{lemma}

\begin{proof}
  This follows from Lemma~\ref{lem:ttau}.
\end{proof}

\subsection{Hyperplane Arrangements}
\label{subsec:hyperplane}

For $s\in\R^r\setminus\{\zero\}$ let 
$H(s):=\narrowsetdef{x\in\R^r}{s^Tx=0}$.  The two connected components
of~$\R^r\setminus H(s)$ are denoted by $H^+(s)$ and $H^-(s)$, where
$s\in H^+(s)$.  For a finite subset $S\subset\R^r\setminus\{\zero\}$ denote
by $ \mathcal{H}(S):=\narrowsetdef{H(s)}{s\in S} $ the \emph{hyperplane
  arrangement} defined by~$S$. The connected components of $
\overline{\mathcal{H}(S)}:=\R^r\setminus\bigcup_{s\in S}H(s) $ are
the \emph{chambers} of $\mathcal{H}(S)$.  We denote the number of
chambers of $\mathcal{H}(S)$ by $\chi(S)$.

\begin{observation}
  Let $C$ be a chamber of $\mathcal{H}(S)$ for some finite subset
  $S\subset\R^r\setminus\{\zero\}$. For each $s\in S$, we have either
  $C\subseteq H^+(s)$ or $C\subseteq H^-(s)$. Defining
  $\varepsilon(C)_s:=+1$ in the first, and $\varepsilon(C)_s:=-1$ in
  the second case, we may assign a \emph{sign vector}
  $\varepsilon(C)\in\{-1,+1\}^S$ to each chamber~$C$ of
  $\mathcal{H}(S)$. This assignment is injective.
\end{observation}

\begin{lemma}
  \label{lem:confarr}
  For each finite subset
  $S\subset\R^r\setminus\{\zero\}$, the
  following equation holds:
  $$
    \#
    \setdef{\varepsilon\in\{-1,+1\}^S}{\zero\not\in\convOp\setdef{\varepsilon_s s}{s\in S}}
    \ =\
  \chi(S)
  $$
\end{lemma}

\begin{proof}
  Let $S\subset\R^r\setminus\{\zero\}$ be finite. By the Farkas-Lemma
  (linear programming duality), for each $\varepsilon\in\{-1,+1\}^S$,
  we have $\zero\not\in\convOp\narrowsetdef{\varepsilon_s s}{s\in S}$
  if and only if there is some $h\in\R^r$ such that $h^T(\varepsilon_s
  s)>0$ holds for all $s\in S$, which in turn is equivalent to
  $$
  h^Ts\ 
  \begin{cases}
    >\ 0 & \text{ if }\varepsilon_s=+1\\
    <\ 0 & \text{ if }\varepsilon_s=-1\\
  \end{cases}
  $$
  for all $s\in S$. Since the latter condition is
  equivalent to $\varepsilon$ being the sign vector of some chamber of
  $\mathcal{H}(S)$, the statement of the lemma follows.
\end{proof}

Lemma~\ref{lem:alphazplus} and Lemma~\ref{lem:confarr}
immediately yield the following result.

\begin{lemma}
  \label{lem:alphachi}
  For $Z_{r,m}^+\in\binom{\vecconf{r}^+}{m}$ chosen uniformly at
  random, we have
  $$
  \alpha(r+1,m)\ =\ \frac{1}{2^m}\cdot\expect{\chi(Z_{r,m}^+)}\ .
  $$
\end{lemma}

The following upper bound on $\chi(\cdot)$ will (via
Lemma~\ref{lem:alphachi}) yield  upper bounds on
$\alpha(\cdot,\cdot)$ that are sufficient for our needs.
We denote $b(p,q):=\sum_{i=0}^p\binom{q}{i}$.

\begin{theorem}[Harding, see Winder {\cite[p.~816]{Win66}}]
  \label{thm:ubchi}
  For $S\in\binom{\R^r\setminus\{\zero\}}{m}$, we have
  $$
  \chi(S)\ \le\ 2b(r-1,m-1)\ .
  $$
\end{theorem}

\subsection{Bounds on \boldmath{$\tau(k,m)$}}
\label{subsec:bounds}


\begin{proposition}
  \label{prop:collect}
  For $0\le m\le 2^k-2$ the following inequality holds:
  $$
    \tau(k,m)\ \le\ \frac{b(k-2,m-1)}{2^{m-1}}
  $$
\end{proposition}

\begin{proof}
  With $r=k-1$, Lemma~\ref{lem:tau2}, Lemma~\ref{lem:alphachi}, and
  Theorem~\ref{thm:ubchi} yield this. 
\end{proof}

In fact, one can prove that, if~$m$ is not too large relative to~$k$,
then the bound of Proposition~\ref{prop:collect} is
asymptotically sharp as~$k$ tends to infinity. Since we do
not need the result here, we omit the proof which (next to the
theorem of Winder's cited in Theorem~\ref{thm:ubchi}) relies on the
fact that the probability of an $l\times l$ matrix with entries from
$\{-1,+1\}$ (chosen uniformly at random) being singular converges to zero
for $l$ tending to infinity (see \cite{KKS95}).

\begin{proposition}
  \label{prop:collect:2}
  For $m(k)\in o\big(2^{\frac{k}{2}}\big)$, we have
  $$
  \lim_{k\rightarrow\infty}\left(\tau(k,m(k))-\frac{b(k-2,m(k)-1)}{2^{m(k)-1}}\right)\ =\  0\ .
  $$
\end{proposition}

\subsection{A Threshold for \boldmath{$\tau(k,m)$}}
\label{subsec:thresh}

For $x\in\R$, let
$$
\Phi(x)\ :=\ \frac{1}{\sqrt{2\pi}}\int_{-\infty}^x e^{-\frac{t^2}{2}}\dd{t}\ ,
$$
i.e., $\Phi$ is the density function of the normal distribution.

\begin{lemma}[de Moivre-Laplace theorem]
  \label{lem:MoivreLaplace}
  For each $\mu\in\R$, the following holds:
  $$
  \lim_{q\rightarrow\infty}
    \frac{b\left(\lfloor\frac{q}{2}+\mu\sqrt{q}\rfloor,q\right)}{2^q}
  \ =\ 
  \Phi(2\mu)
  $$
\end{lemma}

\begin{theorem}
  \label{thm:tauthrzero}
  For each $\varepsilon>0$, we have
  $$
  \lim_{k\rightarrow\infty}\tau\big(k,\lceil(2+\varepsilon)k\rceil\big)\ =\ 0\ .
  $$
\end{theorem}

\begin{proof}
  Let $\varepsilon>0$ be fixed, and define, for each $k$,
  $
  m_{\varepsilon}^+(k):=\lceil(2+\varepsilon)k\rceil
  $.
  
  Let $\delta>0$ be arbitrarily small, and choose $\mu<0$ such that
  \begin{equation}
    \label{eq:tauthr:1}
    \Phi(2\mu)\ <\ \frac{\delta}{2}\ .     
  \end{equation}
  Due to
  $\displaystyle\lim_{k\rightarrow\infty}\frac{m_{\varepsilon}^+(k)}{k}=2+\varepsilon$, we have, for large enough~$k$, 
  \begin{equation}
    \label{eq:tauthr:2}
    k-2\ \le\
    \frac{m_{\varepsilon}^+(k)-1}{2}+\mu\sqrt{m_{\varepsilon}^+(k)-1}\ .
  \end{equation}
  Due to Proposition~\ref{prop:collect}, we have
  \begin{equation}
    \label{eq:tauthr:3}
    \tau(k,m_{\varepsilon}^+(k))
    \ \le\
    \frac{b(k-2,m_{\varepsilon}^+(k)-1)}{2^{m_{\varepsilon}^+(k)-1}}\ .
  \end{equation}
  Since $b(\cdot,\cdot)$ is monotonically increasing in the first
  component, (\ref{eq:tauthr:2}) yields that the right-hand side
  of~(\ref{eq:tauthr:3}) is bounded from above by
  \begin{equation}
    \label{eq:tauthr:4}
    \frac{b\Big(\frac{m_{\varepsilon}^+(k)-1}{2}+\mu\sqrt{m_{\varepsilon}^+(k)-1}\ ,\
            m_{\varepsilon}^+(k)-1\Big)}{2^{m_{\varepsilon}^+(k)-1}}\ .    
  \end{equation}
  By Lemma~\ref{lem:MoivreLaplace} (with~$q$ substituted by $m_{\varepsilon}^+(k)-1$),
  (\ref{eq:tauthr:4}) may be bounded from above
  by $\Phi(2\mu)+\frac{\delta}{2}$ for all large enough~$k$ (because
  of $\displaystyle\lim_{k\rightarrow\infty}m_{\varepsilon}^+(k)=\infty$).
  Thus, from (\ref{eq:tauthr:1}) we obtain
  $$
  \tau(k,m_{\varepsilon}^+(k))
  \ <\
  \delta 
  $$
  for all large enough~$k$.
\end{proof}

Exploiting Proposition~\ref{prop:collect:2}, one can also prove the
following result. It complements Theorem~\ref{thm:tauthrzero}, but
since we will not need it in our treatment, we do not give a proof here.

\begin{theorem}
  \label{thm:tauthrone}
  For each $\varepsilon>0$ we have
  $$
  \lim_{k\rightarrow\infty}\tau\big(k,\lfloor(2-\varepsilon)k\rfloor\big)\ =\ 1
  $$  
\end{theorem}

\section{The Edge Probability \boldmath{$\pi(d,n)$}}
\label{sec:p}

Throughout this section, let the set $X_{d,n}\in\binom{\pmcube{d}}{n}$
be drawn uniformly at random, $P_{d,n}:=\convOp X_{d,n}$, and let
$\{v,w\}\in\binom{X_{d,n}}{2}$ be chosen uniformly at random as well.
Our aim is to determine the probability
$$
\pi(d,n)\ :=\ \prob{\{v,w\}\in\edges{P_{d,n}}}\ .
$$
Let us further denote
$$
\pi_k(d,n)\ :=\ \probcond{\{v,w\}\in\edges{P_{d,n}}}{\dist{v}{w}=k}\ .
$$
Since $\{v,w\}$ is uniformly distributed over $\binom{\pmcube{d}}{2}$,
the distance $\dist{v}{w}$ has the same distribution as the number of
positive components of a point chosen uniformly at random from
$\pmcube{d}\setminus\{-\one\}$.  Therefore, the following equation
holds.

\begin{lemma}
  \label{lem:pk}
  $$
    \pi(d,n)\ =\ \frac{1}{2^d-1}\sum_{k=1}^d \binom{d}{k}\pi_k(d,n)
  $$
\end{lemma}

The following result, stating that $\pi(d,\cdot)$ is monotonically
increasing, is quite plausible. Its straightforward proof is omitted
here.

\begin{proposition}
  \label{prop:pmonotone}
  The function $\pi(d,\cdot)$ is monotonically decreasing, i.e., 
  for $3\le n\le 2^d-1$, we have
  $
  \pi(d,n)>\pi(d,n+1)
  $.
\end{proposition}

The next result implies part~(i) of Theorem~\ref{thm:dens} (see
the remarks at the end of Section~\ref{sec:intro}).

\begin{proposition}
  \label{prop:pthrlow}
  For each $\varepsilon>0$, we have
  $$
  \lim_{d\rightarrow\infty}\ \pi\left(d,\left\lfloor 2^{\left(\frac{1}{2}-\varepsilon\right)d}\right\rfloor\right)
  \ =\ 1\ .
  $$
\end{proposition}

\begin{proof}
  Let $\varepsilon>0$, and define
  $
    n^-_{\varepsilon}(d):=\big\lfloor 2^{\left(\frac{1}{2}-\varepsilon\right)d}\big\rfloor
  $.  For
  each $\mu>0$, denote
  $$
  K^{\le}_{\mu}(d)\ :=\ \narrowsetdef{k\in\Z}{1\le k\le \frac{d}{2}+\mu\sqrt{d}}
  $$
  and 
  $$
  \pi^-_{\mu}(d)\ :=\ \min\bigsetdef{\pi_k(d,n^-_{\varepsilon}(d))}{k\in K^{\le}_{\mu}(d)}\ .
  $$
  Then, due to Lemma~\ref{lem:pk}, we have
  $$
    \pi(d,n^-_{\varepsilon}(d))\ \ge\ 
    \sum_{k\in K^{\le}_{\mu}(d)}\frac{\binom{d}{k}}{2^d}\cdot \pi^-_{\mu}(d)\ .
  $$
  For every $\nu>0$, this implies (by Lemma~\ref{lem:MoivreLaplace})
  that
  \begin{equation}
    \label{eq:pthrlow:05}
    \pi(d,n^-_{\varepsilon}(d))\ \ge\ (\Phi(2\mu)-\nu)\cdot \pi^-_{\mu}(d)
  \end{equation}
  holds for all large enough~$d$. Therefore, it remains to prove, for
  all $\mu>0$,
  \begin{equation}
    \label{eq:pthrlow:1}
    \lim_{d\rightarrow\infty}\pi^-_{\mu}(d)\ =\ 1\ .
  \end{equation}
  With
  $$
  \xi_k\ :=\ \probcond{X_{d,n^-_{\varepsilon}(d)}\cap\pmcubestarOp(v,w)=\varnothing}{\dist{v}{w}=k}\ ,
  $$
  we have, for each $k\in K^{\le}_{\mu}(d)$, 
  \begin{equation}
    \label{eq:pthrlow:2}
    \pi_k(d,n^-_{\varepsilon}(d))\ \ge\ \xi_k\ \ge\
    \xi_{\left\lfloor\frac{d}{2}+\mu\sqrt{d}\right\rfloor}
  \end{equation}
  (see Lemma~\ref{lem:adjacency}).
  Clearly,
  $$
    \expectcond{\#(X_{d,n^-_{\varepsilon}(d)}\cap\pmcubestarOp(v,w))}{\dist{v}{w}=k}
    \ =\
    \frac{2^k-2}{2^d-2}\cdot(n^-_{\varepsilon}(d)-2)\ ,
  $$
  and thus, the estimation
  $$
    \expectcond{\#(X_{d,n^-_{\varepsilon}(d)}\cap\pmcubestarOp(v,w))}{\dist{v}{w}=k}
    \ \le\
    2^{k-(\frac{1}{2}+\varepsilon)d}\ ,
  $$
  hold for each $k$. By Markov's inequality, this
  implies
  \begin{equation}
    \label{eq:pthrlow:4}
    \probcond{\#(X_{d,n^-_{\varepsilon}(d)}\cap\pmcubestarOp(v,w))\ge d\cdot2^{k-(\frac{1}{2}+\varepsilon)d} }{\dist{v}{w}=k}
    \ \le\ \frac{1}{d}
  \end{equation}
  for each~$d$ and~$k$. For
  $k=\big\lfloor\frac{d}{2}+\mu\sqrt{d}\big\rfloor$,
  (\ref{eq:pthrlow:4}) yields
  \begin{multline}
    \label{eq:pthrlow:5}
    \probcond{\#(X_{d,n^-_{\varepsilon}(d)}\cap\pmcubestarOp(v,w))\ge d\cdot2^{\mu\sqrt{d}-\varepsilon d} }%
             {\dist{v}{w}=\big\lfloor\frac{d}{2}+\mu\sqrt{d}\big\rfloor}\\
    \le\ \frac{1}{d}    
  \end{multline}
  for all~$d$. Since $d\cdot2^{\mu\sqrt{d}-\varepsilon d}< 1$
  holds for large enough~$d$, (\ref{eq:pthrlow:5}) implies
  $
  \xi_{\left\lfloor\frac{d}{2}+\mu\sqrt{d}\right\rfloor}\ge 1-\frac{1}{d}
  $
  for large enough~$d$.
  Therefore, 
  $$
  \lim_{d\rightarrow\infty}\xi_{\left\lfloor\frac{d}{2}+\mu\sqrt{d}\right\rfloor}\ =\ 1
  $$
  holds, which, by (\ref{eq:pthrlow:2}), finally implies (\ref{eq:pthrlow:1}).
\end{proof}

The next result yields part~(ii) of Theorem~\ref{thm:dens} (see
the remarks at the end of Section~\ref{sec:intro}).

\begin{proposition}
  \label{prop:pthrup}
  For each $\varepsilon>0$, we have
  $$
  \lim_{d\rightarrow\infty}\pi\left(d,\left\lceil 2^{\left(\frac{1}{2}+\varepsilon\right)d}\right\rceil\right)
  \ =\ 0\ .
  $$
\end{proposition}

\begin{proof}
  Let $\varepsilon>0$, and define
  $
    n^+_{\varepsilon}(d):=\left\lceil 2^{\left(\frac{1}{2}+\varepsilon\right)d}\right\rceil
  $.
  For each $\mu>0$, denote
  $$
  K^{\ge}_{\mu}(d)\ :=\ \narrowsetdef{k\in\Z}{\frac{d}{2}-\mu\sqrt{d}\le
    k\le d}\ ,
  $$
  and define 
  \begin{equation}
    \label{eq:pthrup:0}
    \pi^+_{\mu}(d)\ :=\ \max\setdef{\pi_k(d,n^+_{\varepsilon}(d))}{k\in K^{\ge}_{\mu}(d)}\ .
  \end{equation}
  Then, due to Lemma~\ref{lem:pk}, we have
  $$
    \pi(d,n^+_{\varepsilon}(d))\ \le\ 
    2\cdot\hspace{-3mm}\sum_{k=1}^{\lfloor
      \frac{d}{2}-\mu\sqrt{d}\rfloor}\frac{\binom{d}{k}}{2^d}\ +\ \pi^+_{\mu}(d)\ .
  $$
  Thus, for every $\nu>0$, by Lemma~\ref{lem:MoivreLaplace},
  $$
    \pi(d,n^+_{\varepsilon}(d))\ \le\ \Phi(-2\mu)+\nu+\pi^+_{\mu}(d)
  $$
  holds for all large enough~$d$. Therefore, it remains to prove, for
  all $\mu>0$,
  \begin{equation}
    \label{eq:pthrup:1}
    \lim_{d\rightarrow\infty}\pi^+_{\mu}(d)\ =\ 0\ .
  \end{equation}
  For $k\in\{1,\dots,d\}$ and
  $m\in\{0,\dots,2^k-2\}$, we define
  $$
  \xi_k(m)\ :=\ \probcond{\#(X_{d,n^+_{\varepsilon}(d)}\cap\pmcubestarOp(v,w))=m}{\dist{v}{w}=k}
  $$
  (i.e., $\xi_k(0)=\xi_k$ in the proof of
  Proposition~\ref{prop:pthrlow}). 
  Then we have (see Lemma~\ref{lem:tau})
  \begin{equation}
    \label{eq:pthrup:25}
    \pi_k(d,n^+_{\varepsilon}(d))\ =\ \sum_{m=0}^{2^k-2}\xi_k(m)\tau(k,m)\ .
  \end{equation}
  Since $\tau(k,\cdot)$ is monotonically non-increasing by
  Lemma~\ref{lem:taumono}, we thus can estimate
  $$
    \pi_k(d,n^+_{\varepsilon}(d))\ \le\ \sum_{m=0}^{3k-1}\xi_k(m)\ +\ \tau(k,3k)\ ,
  $$
  for each $k\in K^{\ge}_{\mu}(d)$.  This yields, again for for each
  $k\in K^{\ge}_{\mu}(d)$,
  \begin{multline}
    \label{eq:pthrup:2}
    \pi_k(d,n^+_{\varepsilon}(d))\ \le\ 3d\cdot\max\setdef{\xi_k(m)}{0\le m\le 3d-1}\ \\
    +\ \max\setdef{\tau(k',3k')}{k'\in K^{\ge}_{\mu}(d)}\ .
  \end{multline}
  According to Theorem~\ref{thm:tauthrzero}, 
  $$
  \lim_{d\rightarrow\infty}\max\narrowsetdef{\tau(k',3k')}{k'\in
    K^{\ge}_{\mu}(d)}\ =\ 0
  $$
  holds. Hence, by (\ref{eq:pthrup:2}) and
  (\ref{eq:pthrup:0}), equation~(\ref{eq:pthrup:1}) can be proved by showing 
  \begin{equation}
    \label{eq:pthrup:3}
    \lim_{d\rightarrow\infty}\ 
      \big(3d\cdot\max
        \setdef{\xi_k(m)}{0\le m\le 3d-1,\ k\in K^{\ge}_{\mu}(d)}
      \big)\
    =\ 0\ .
  \end{equation}  
  Let us first calculate (using the notation $\lowfac{a}{b}:=a(a-1)\cdots(a-b+1)$)
  \begin{eqnarray}
    \label{eq:pthrup:4}
    \xi_k(m)\ &=\ & {\displaystyle
                      \frac{\binom{2^k-2}{m}\binom{2^d-2^k}{n^+_{\varepsilon}(d)-m-2}}%
                           {\binom{2^d-2}{n^+_{\varepsilon}(d)-2}}
                    }\nonumber\\
              &=  & {\displaystyle
                           \binom{2^k-2}{m}
                      \cdot\frac{\lowfac{2^d-2^k}{n^+_{\varepsilon}(d)-m-2}}%
                                {\lowfac{2^d-2}{n^+_{\varepsilon}(d)-2}}
                      \cdot\frac{(n^+_{\varepsilon}(d)-2)!}{(n^+_{\varepsilon}(d)-m-2)!}
                    }\ ,  
  \end{eqnarray}
  where the left, the middle, and the right factor of (\ref{eq:pthrup:4}) may be
  bounded from above by
  $
  (2^d)^m
  $,
  $
         (2^d)^2\cdot\left(\frac{2^d-2^k}{2^d}\right)^{n^+_{\varepsilon}(d)}
  $,
  and
  $
  (2^d)^{m}
  $,
  respectively. Thus, we obtain, for $0\le m\le 3d-1$,
  \begin{equation}
    \label{eq:pthrup:5}
    \xi_k(m)\ \le\ 2^{\const\cdot d^2}\cdot\left(1-\frac{1}{2^{d-k}}\right)^{n^+_{\varepsilon}(d)}\ .
  \end{equation}
  For $k\in K^{\ge}_{\mu}(d)$, we have
  \begin{eqnarray}
    \label{eq:pthrup:6}
    \left(1-\frac{1}{2^{d-k}}\right)^{n^+_{\varepsilon}(d)}
      \ &\le\ & \left(1-\frac{1}{2^{\frac{d}{2}+\mu\sqrt{d}}}\right)%
                  ^{2^{\left(\frac{1}{2}+\varepsilon\right)d}}\nonumber\\
        &=    & \left[%
                  \left(1-\frac{1}{2^{\frac{d}{2}+\mu\sqrt{d}}}\right)%
                    ^{2^{\frac{d}{2}+\mu\sqrt{d}}}
                \right]%
                  ^{2^{\varepsilon d-\mu\sqrt{d}}}\ .
  \end{eqnarray}
  For~$d$ tending to infinity, the expression in the square brackets
  of (\ref{eq:pthrup:6}) converges to $\frac{1}{\Euler}<\frac{1}{2}$
  (where~$\Euler=2.7182\cdots$ is Euler's constant). Therefore,
  (\ref{eq:pthrup:6}) and (\ref{eq:pthrup:5}) imply
  $
    \xi_k(m)\le 2^{\const\cdot d^2}\cdot(1/2)^{2^{\varepsilon d-\mu\sqrt{d}}}
  $
  (for $k\in K^{\ge}_{\mu}(d)$, $0\le m\le 3d-1$, and for large
  enough~$d$). This finally yields (\ref{eq:pthrup:3}), and therefore completes the proof.
\end{proof}

\section{Remarks}
\label{sec:rem}

The threshold for the function $\tau(\cdot,\cdot)$
described in Theorems~\ref{thm:tauthrzero} and~\ref{thm:tauthrone}
is much sharper than we needed for our purposes (proof of
Proposition~\ref{prop:pthrup}). The sharper result may, however, be
useful in investigations of more structural properties of the graphs
of random 0/1-polytopes. A particularly interesting such question is
whether these graphs have good expansion properties with high
probability.

\section*{Acknowledgements}
We thank one of the referees for several suggestions that helped to
improve the presentation.


\end{document}